\documentclass[11pt]{amsproc}
\usepackage{amstext,amsmath,amssymb,amsthm, esint}
\usepackage{latexsym}
\usepackage{exscale}
\usepackage{tikz}

\usepackage{color}

 \RequirePackage{srcltx}

\numberwithin{equation}{section}
\renewcommand\r{\rangle}
\renewcommand\l{\langle}

\newcommand\dsize{\displaystyle}

\newcommand\R{\mathbb{R}}

\newcommand\Z{\mathbb{Z}}
\newcommand\N{\mathbb{N}}

\newcommand\G{{\bf\Gamma}}

\newcommand\B{\mathcal{B}}

\newtheorem{Thm}{Theorem}[section]
\newtheorem{Lemma}[Thm]{Lemma}

\theoremstyle{remark}

\begin{document}

\title[]{    Methods of Construction of Exponential Bases on Planar Domains Revised}
\author{Oleg Asipchuk}
\address{Oleg Asipchuk,  University of Cincinnati,
	Department of Mathematical Sciences,
	Cincinnati, OH 45221, USA}
\email{asipchah@ucmail.uc.edu}
\author{Laura De Carli}
\address{L.~De Carli, Florida International University,
Department of Mathematics,
Miami, FL 33199, USA}
\email{decarlil@fiu.edu}
	\maketitle
 \begin{abstract}
 In \cite{AD}, we surveyed classical and recent results on the construction of exponential bases on planar domains with explicit control over the associated frame bounds. We also constructed well-conditioned exponential bases on certain planar sets that multi-tile the plane. In this revision, we identify an error in the choice of basis used in the proof of Theorem 3.11 and we provide a new proof using a simpler exponential basis.  
\end{abstract} 

\section{Introduction}

\begin{figure}[ht] 
	\begin{center}
		\begin{tikzpicture}
			\draw[black, ->] (0, 0)--(5, 0);
			\draw[black, ->] (0, 0)--(0, 5);
			\draw [black ] (0,0) --(4,0)--(4,1)--(0,1)--(0,0);
			\fill[lightgray ](0,0) --(4,0)--(4,1)--(0,1)--(0,0);
			\draw [black] (0,1) --(0,2)--(2,4)--(4,2)--(4,1)--(2,3)--(0,1);
			
			\fill[lightgray ](0,1) --(0,2)--(2,4)--(4,2)--(4,1)--(2,3)--(0,1);
			\draw    (-0.5,-0.5) node [black]{{\small $(0,0)$}};
			\draw    (4.5,-0.5) node [black]{{\small $(1, 0)$}};
			\draw    (2,3.5) node [black]{{\small $R_2$}};
			\draw    (2,0.5) node [black]{{\small $R_1$}};
			\draw    (-0.5,1) node [black]{{\small $(0,\frac{1}{4})$}};
			\draw (-0.5, 2)  node [black]{{\small $(0,\frac{1}{2})$}};
			\draw    (4.5,1) node [black]{{\small $(1, \frac 14)$}};
			\draw    (2,4.5) node [black]{{\small $(\frac{1}{2},1)$}};
			\draw    (2,2.3) node [black]{{\small $(\frac{1}{2},\frac{3}{4})$}};
			
		\end{tikzpicture}\caption { The domain $D=R_1\cup R_2$ }
	\end{center}
	\label{f1}
\end{figure}
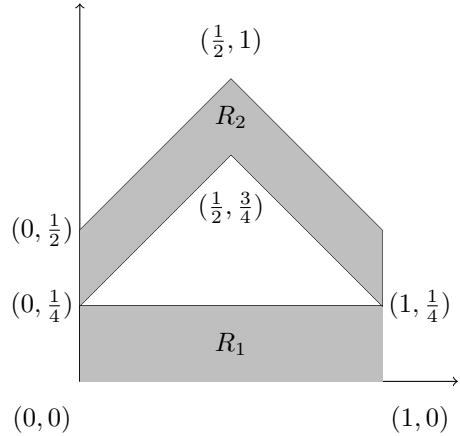

One of the main results in our  paper \cite{AD}   is the construction of an explicit   well-conditioned exponential basis on the domain  in figure \ref{f1}.

In \cite[Theorem 3.11]{AD},  we proved that  the set $\B=\bigcup_{N\in\N}\B'_{N} $, where
$$
\B'_N =\!\!\!\bigcup_{j=0}^{2^{2(N-1)}-1}\!\!\left\{e^{2 \pi i ((2^Nn+\frac{4 j}{2^N})x+(2^Nm+\frac{4 j}{2^N})y)}, \  e^{2 \pi i ((2^Nn+\frac{4 j}{2^N}+\frac{2}{2^N})x+(2^N m+\frac{4 j}{2^N}+\frac{2}{2^N})y)}\right\}_{n, m\in\Z}
$$
is a Riesz basis for $L^2(D )$, with frame  bounds $ 
B=\frac{2+\sqrt{2} }{4}$   {and } $
A= \frac{2-\sqrt{2} }{4}. 
$ 

 \medskip
  The existence of an exponential basis on $L^2(D)$ was not the novelty of our theorem. Indeed, it is not too difficult to verify that there exists a lattice   $L\subset \R^2$ for which the collection $\{D + \ell : \ell \in L\}$ covers  the plane  exactly  two times,  up to  a set of measure zero, and every set that multi-tiles the plane   has an exponential basis.      However, the available general results do not necessarily provide an explicit construction of such a basis or estimates for its Riesz bounds. The methods developed in \cite{AD}, when applicable, provide explicit exponential bases together with quantitative bounds.
  
  Unfortunately,  the exponential  set $\B$  that we used in the proof of Theorem 3.11     contains  frequencies $$ \{(2^N n+\frac{4j}{2^N}, \ 2^N m + \frac{4j}{2^N} )\}_{n, m\in\Z, N \ge 2 \atop{0\leq j\leq 2^{N-2}-1}},$$  and   when $n=m =0$ and $j=1$, the sequence  $\{ ( \frac{4 }{2^N}  \   \frac{4 }{2^N} )\}_{ {N \ge 2}}$ converges to $(0,0)$. 
  Uniform separation of frequencies is
 a necessary condition for a set of exponentials to form a Riesz basis, and so our set  $\B$ cannot be a  Riesz basis.  
   
   
   \medskip
   We  should also  mention   a typo in  the  expression for $B$ in  the statement of Theorem 2.8.   We should have $B\leq 16$, not $B\leq 14$.

   \medskip
   The correct version of  \cite[Theorem 3.11]{AD} is the following

 \begin{Thm}\label{main-theorem}
 	Let \(D\) be the domain in Figure \ref{f1}. The exponential
 	system
 	\begin{equation}\label{main-spectrum}
 	 	\B  =\!\!\! \left\{e^{2 \pi i ( nx+4my)   }\right\}_{n, m\in\Z}\cup \left\{e^{2 \pi i ( nx+ (4m+1)y   }\right\}_{n, m\in\Z}\
 	 \end{equation}
 	is a Riesz basis for \(L^2(D)\), with Riesz bounds  
 	\begin{equation}\label{main-bounds}
 		A=\frac{2-\sqrt2}{4},
 		\qquad
 		B=\frac{2+\sqrt2}{4}.
 	\end{equation}
 \end{Thm}

Letting $E(\Lambda):=\{e^{2\pi i (\lambda_1   x + \lambda_2 y)}\}_{(\lambda_1, \lambda_2)\in\Lambda }$, we can also  write   $$\mathcal B
=
E\bigl(\mathbb Z\times4\mathbb Z\bigr)
\cup
E\bigl(\mathbb Z\times(4\mathbb Z+1)\bigr).$$

As in the original paper, our proof relies on Lemmas \ref{L-limit-bases} and  Theorem~\ref{T1-N-interv} stated below, but  now we use   the fixed exponential system $\B$  defined  in \eqref{main-spectrum},  which is a Riesz basis for every domain $D_N$  that approximate $D$.
 
\section{ Constructing an exponential basis on $D$}  
The proof of Theorem \ref{main-theorem} is based on the  following results, which are  stated and proved in \cite{AD}.

 \begin{Lemma}\label{L-limit-bases}
 	Let  $\Omega, \,D \subset \R^d$ be  bounded and measurable, with $\Omega \supset D$. Let    $\{D_n\}_{n=1}^\infty $   be  a  family of measurable subsets of $\Omega$, with   $D_n\cap D\uparrow D$ and 
 	$D_n\cup D\downarrow D$.

 	Let   $\B_n \subset L^2(\Omega)\cap L^\infty(\Omega) $ be a Riesz basis of $ L^2(D_n  ) $ with frame  bounds   $A_n$  and $B_n$, and  let $\B=\cup_{n\in\N} \B_n$. Assume   $\ \B_n\subset\B_{n+1}$  for every $n\in\N$ and that 
 	$$ 0<  \alpha\leq A_n\leq B_n\leq \beta <\infty, $$ where   $\alpha,\beta> 0$ are independent of $n$.
 	Assume also that  the $L^\infty$ norm of the functions in  $\B $  is  bounded above by a constant independent of $n$. 
 	Then, $\B $ is a Riesz basis of $L^2(D)$ with frame  bounds  
 	$  \alpha $ and $   \beta.$ 
 \end{Lemma}
 
We denote     $Q_d  $  the  unit cube $[ -\frac 12,  \,\frac 12)^d$; given  $p_1,\, ...,\, p_N\in \R^d $ such that 
$Q_d+  p_k\cap Q_d+  p_h=\emptyset$ when $h\ne k$, we let 
$$ T(  p_1, \, ...,  p_N):=\dsize \bigcup_{k=1}^N (Q_d+  p_k).$$

 \begin{Thm}\label{T1-N-interv}
 With the definition and the notation previously introduced, the set   
 	$$ \B= \B( \delta_1, \, ...,\,   \delta_N):=  \bigcup_{ j=1} ^N \{e^{ 2\pi i \l   n+  \delta_j,\,  x\r}\}_{  n\in\Z^d} $$
 	is a Riesz basis of  $L^2(T(p_1,\, ...,\, p_N))$  if and only if  the matrix
 	\begin{equation}\label{def-Gamma}{\bf \Gamma}=\{ e^{2\pi i \l  \delta_j, p_k\r}\}_{1\leq j,k\leq N} \end{equation} is nonsingular.  
 	The optimal frame bounds of $\B$  are the maximum and minimum singular values of $\G$.
 \end{Thm}

  \subsection{ Proof of Theorem \ref{main-theorem}}
 
   Recall that we have denoted \mbox{$D=R_1\cup R_2$}, where $R_1$ is the rectangle at the bottom of Figure \ref{f2}. and $R_2=D-R_1$. 
   Our proof  uses most of the results listed above and consists of several steps:
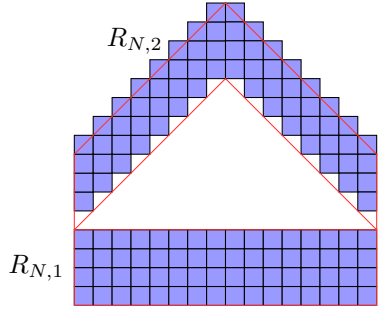
\begin{figure}[ht]
	\begin{center}
		\begin{tikzpicture}

			\foreach \y in {0,1,...,3}{
				\foreach \x in {0,1,...,15} {
					\fill[blue!40!white,draw=black] (\x/4,\y/4) rectangle (\x/4+1/4,\y/4+1/4);
			}}

			\foreach \x in {0,1,...,7}{
				\foreach \y in {0,1,...,3} {
					\fill[blue!40!white,draw=black] (\x/4,\y/4+1+1/4+\x/4) rectangle (\x/4+1/4,\y/4+1+1/4+\x/4+1/4);
			}}
			
			\foreach \x in {0,1,...,7}{
				\foreach \y in {0,1,2,...,3} {
					\fill[blue!40!white,draw=black] (\x/4+2,4-1/4-\y/4-\x/4) rectangle (\x/4+2+1/4,4-1/4-\y/4-\x/4+1/4);
			}}

			\draw    (0.8,3.5) node [black]{{\small $R_{N,2} $}};
			\draw    (-0.5,0.5) node [black]{{\small$ R_{N,1} $}};
			
			\draw [red!80!white] (0,0) --(4,0)--(4,1)--(0,1)--(0,0);
			\draw [red!80!white] (0,1) --(0,2)--(2,4)--(4,2)--(4,1)--(2,3)--(0,1);  
		\end{tikzpicture}
	\end{center}
	\caption{Approximation of $D$ by $D_N$}\label{f2}
\end{figure}

\medskip
\noindent
{\textbf{a) Approximation}}: We consider the  multi-rectangles $\{D_{ N}\}_{N=2}^\infty = $ \newline$\{R_{N,1} \cup R_{N,2} \}_{N=2}^\infty $   with sides $\frac{1}{2^N}$  and lower-left vertices  in the sets $A_{N,1} \cup A_{N,2} $, where 
$ 
A_{N,1}  = \left\{\left(\frac{k}{2^N},\frac{l}{2^N}\right)\right\}_%
{\genfrac{}{}{0pt}{2}{0\leq k \leq 2^N-1}{0\leq l \leq \frac{2^N}{4}-1}}
$ 
and 
$$
 A_{N,2}  = \left\{\left(\frac{k}{2^N},\, \frac{1}{4}+\frac{k+l+1}{2^N}\right),\ 
 \left(\frac{1}{2}+\frac{k}{2^N},\frac{2^N-2-k-l}{2^N}
 \right) \right\}_{\genfrac{}{}{0pt}{2}{0\leq k \leq \frac{2^N}{2}-1}{0\leq l \leq \frac{2^N}{4}-1}}.
 $$
See Figure \ref{f2}.
 
 For every $N\ge 2$, the  union of squares in the set $A_{N,1} $   forms the rectangle $ R_{1, N}$ at the bottom of Figure \ref{f2}. Since $R_{N,1} =R_1=R_{N+1,1} $
 for every $N\ge 2$, no approximation is required for the lower rectangle.

 We constructed the dyadic multi-rectangles $R_{N, 2}$,  in such a way that the portion of $D_N$ lying inside $D$ increases with $N$,  while the portion of $D_N$ lying outside $D$ decreases. More precisely,
$$
 D_N\cap D\subset D_{N+1}\cap D,
 \qquad
  D_{N+1}\cup D\subset D_N\cup D.
$$
 Furthermore, the width and height of the boundary squares are $2^{-N}$, and so  the boundaries of the  "staircases" $R_{N,2}$  converge to the sloped sides of $D$.  Consequently,
$$\bigcup_{N=2}^{\infty}(D_N\cap D)=D
  \quad\text{and}\quad
  \bigcap_{N=2}^{\infty}(D_N\cup D)=D,
$$
 up to sets of measure zero. Moreover,  our definitions  ensure  that 
 \mbox{$D\cap R_{N,2} \subset D\cap R_{ N+1, 2}$.}
 Therefore,
$$
  D_N\cap D\uparrow D,\qquad D_N\cup D\downarrow D,
$$
  as required in Lemma \ref{L-limit-bases}.  

\medskip
\noindent
\textbf{b) Dilation}: Next, we apply  the linear transformation  
$L(x,y)=(2^Nx,2^Ny)$  to obtain unions of unit squares
 $\{	\tilde{D}_{N}\}_{N=2}^\infty = \{	\tilde{R}_{N,1} \cup 	\tilde {R}_{N,2} \}_{N=2}^\infty $  whose lower-left vertices  in   $\tilde {A}_{N,1} \cup 	\tilde {A}_{N,2} $, where
$  
	\tilde {A}_{N,1}  = \left\{\left(k,l\right)\right\}_{\genfrac{}{}{0pt}{2}{0\leq k \leq 2^N-1}{0\leq l \leq 2^{N-2}-1}}
$ 
and 
$$
	\tilde {A}_{N,2} \! =\! \left\{\left(k,2^{N-2}+k+l+1\right) , \ \left(2^{N-1}+k,2^N-2-k-l\right)\right\}_{\genfrac{}{}{0pt}{2}{0\leq k \leq 2^{N-1}-1}{0\leq l \leq 2^{N-2}-1}}.
$$

 We use Theorem \ref{T1-N-interv} to prove that the exponential system 
\begin{equation}\label{E-systemN}
 \mathcal B_N
 =
 E\bigl(\frac{\mathbb Z}{2^N}\times \frac{4\mathbb Z}{2^{N}}\bigr)
 \cup
  E\bigl(\frac{\mathbb Z}{2^N}\times \frac{4 \mathbb Z+1}{2^{N}}\bigr)
\end{equation}
is a Riesz basis for $L^2(\tilde{D}_{N})$. To prove this, we estimate the singular values of the matrix
$$
{\bf \Gamma}=\left(e^{2\pi i\langle \delta_j,p_k\rangle}\right)_{1\leq j,k\leq 2^{2N-1}},
$$
where \(p_k\) are the elements of \(\widetilde{A}_{N,1}\cup \widetilde{A}_{N,2}\), and
\(\delta_j\) are the elements of 
$$\Delta_N=\Delta_{1,N}\bigcup \Delta_{2,N}=\left\{\left( \frac{k}{2^N},\frac{4 l}{2^N}\right)\right\}_{\genfrac{}{}{0pt}{2}{0\leq k \leq 2^N-1}{0\leq l \leq 2^{N-2}-1}}\bigcup \left\{\left( \frac{k}{2^N},\frac{4 l+1}{2^N}\right)\right\}_{\genfrac{}{}{0pt}{2}{0\leq k \leq 2^N-1}{0\leq l \leq 2^{N-2}-1}}.$$
  
\textbf{c) Estimating the singular values}: We now estimate the singular values of $\Gamma$ by analyzing its Gram matrix. Consider two columns $\vec c_1(p_1,p_2)$ and $\vec c_2(p'_1,p'_2)$ of $\Gamma$ corresponding to distinct points
$  
(p_1,p_2),(p'_1,p'_2)\in \widetilde A_{N,1}\cup\widetilde A_{N,2}.
$ 
The inner product of these two vectors is
\begin{align*}
	\l \vec c_1,\vec c_2\r
	&=
	\sum_{k=0}^{2^N-1}\sum_{l=0}^{2^{N-2}-1}
	\left(
	e^{2\pi i\left((p_1-p'_1)\frac{k}{2^N}
		+(p_2-p'_2)\frac{4l}{2^N}\right)}
	+
	e^{2\pi i\left((p_1-p'_1)\frac{k}{2^N}
		+(p_2-p'_2)\frac{4l+1}{2^N}\right)}
	\right)
	\\
	&=
	\left(
	\sum_{k=0}^{2^N-1}
	e^{2\pi i(p_1-p'_1)\frac{k}{2^N}}
	\right)
	\left(
	1+e^{2\pi i\frac{p_2-p'_2}{2^N}}
	\right)
	\sum_{l=0}^{2^{N-2}-1}
	e^{2\pi i(p_2-p'_2)\frac{4l}{2^N}}.
\end{align*}

Since the first coordinates of the points of
$\widetilde A_{N,1}\cup\widetilde A_{N,2}$ belong to
$\{0,\ldots,2^N-1\}$, the first sum vanishes unless
$
p_1=p'_1.
$
If $p_1=p'_1$, then
\begin{align*}
	\l \vec c_1,\vec c_2\r
	&=
	2^N
	\left(
	1+e^{2\pi i\frac{p_2-p'_2}{2^N}}
	\right)
	\sum_{l=0}^{2^{N-2}-1}
	e^{2\pi i(p_2-p'_2)\frac{4l}{2^N}}.
\end{align*}
The last sum vanishes unless
$  
p_2-p'_2\equiv 0\pmod{2^{N-2}}.
$
Thus
$ 
\l \vec c_1,\vec c_2\r=0
 $
unless
$
p_1=p'_1$ and 
$
p_2-p'_2\equiv0\pmod{2^{N-2}}.
$
Equivalently, two distinct columns can have a nonzero inner product only when the corresponding points are congruent modulo the lattice
$$
\mathcal L_N:=2^N\mathbb Z\times2^{N-2}\mathbb Z.
$$
Each coset of $\mathcal L_N$ contains exactly one point of
$\widetilde A_{N,1}$ and exactly one point of
$\widetilde A_{N,2}$. Indeed, $\widetilde A_{N,1}$ is a complete
set of representatives modulo $\mathcal L_N$. For
$\widetilde A_{N,2}$, the first family of points provides exactly
one representative for every coset whose first coordinate modulo
$2^N$ belongs to
$
\{0,\ldots,2^{N-1}-1\},
$
while the second family gives exactly one representative for each of the remaining cosets.
 
Consequently, the columns of $\Gamma$ can be partitioned into pairs,
where two columns belong to the same pair  when the
corresponding points are congruent modulo $\mathcal L_N$.

The previous computation shows that columns belonging to distinct
pairs are orthogonal. Hence, after a suitable permutation of the
columns, the Gram matrix $\Gamma^*\Gamma$ is block diagonal, with
$2\times2$ blocks corresponding to these pairs. Since the
characteristic polynomial of a block-diagonal matrix is the product
of the characteristic polynomials of its diagonal blocks, the
eigenvalues of $\Gamma^*\Gamma$ are  the eigenvalues of
these $2\times2$ blocks, counted with multiplicity.

Consider one such pair and write
$$
\Gamma_{(p,p')}
:=
\big[
\vec c(p_1,p_2); \,
\vec c(p'_1,p'_2)
\big].
$$
Since $p$ and $p'$ are distinct and congruent modulo $\mathcal L_N$,
we have
$$
p_1=p'_1,
\qquad
p_2-p'_2=a\,2^{N-2}.
$$
Interchanging $p$ and $p'$ if necessary, we may assume that
$
a\in\{1,2,3\}.
$
It follows from the preceding computation that the corresponding
$2\times2$ block of the Gram matrix is
$$
\Gamma_{(p,p')}^*\Gamma_{(p,p')}
=
2^{2N-2}
\begin{pmatrix}
	2 & 1+e^{2\pi ia/4}\\
	1+e^{-2\pi ia/4} & 2
\end{pmatrix}.
$$
Its eigenvalues are
$$
2^{2N-2}
\left(
2\pm\left|1+e^{2\pi ia/4}\right|
\right).
$$
Since $a\in{1,2,3}$, we have that
$
\left|1+e^{2\pi ia/4}\right|\leq\sqrt2,  
$
and so
$$
2^{2N}\frac{2-\sqrt2}{4}
\leq
\sigma^2
\leq
2^{2N}\frac{2+\sqrt2}{4}
$$
for every singular value $\sigma$ of $\Gamma$. In particular,
$$
\sigma_{\min}^2(\Gamma)
\geq
2^{2N}\frac{2-\sqrt2}{4},
\qquad
\sigma_{\max}^2(\Gamma)
\leq
2^{2N}\frac{2+\sqrt2}{4}.
$$

Thus, by Theorem~2.2, $\mathcal B_N$ is a Riesz basis for
$L^2(\widetilde D_N)$ with Riesz bounds
$$
\widetilde A_N
\geq
2^{2N}\frac{2-\sqrt2}{4},
\qquad
\widetilde B_N
\leq
2^{2N}\frac{2+\sqrt2}{4}.
$$

Finally, contracting by the factor $2^{-N}$, the system $\mathcal B_N$
in \eqref{E-systemN} becomes the fixed system $\mathcal B$ in
\eqref{main-spectrum}, while its Riesz bounds are multiplied by
$2^{-2N}$. Hence, $\mathcal B$ is a Riesz basis for $L^2(D_N)$,
for every $N\geq2$, with the uniform bounds
$$
A=\frac{2-\sqrt2}{4},
\qquad
B=\frac{2+\sqrt2}{4}.
$$
The conclusion now follows from Lemma~2.1 and Step~a).

  \medskip

  ACKNOWLEDGMENTS:  
   We are indebted to Alex Iosevich for pointing out the flaws in our original manuscript.

  \end{document}